\documentclass[12pt]{article}%[10pt]{article}%\documentclass
\usepackage{amsmath,latexsym,epsf}
\usepackage{amsfonts}
\usepackage{amssymb}
\usepackage{mathrsfs} %no en windows
\usepackage{amsmath}
\usepackage[english]{babel} %si en windows
\usepackage[latin1]{inputenc}
\usepackage[all]{xy}
\begin{document}

%% ----------------------------------------------------------------------
\newcommand{\nc}{\newcommand}
\def\PP#1#2#3{{\mathrm{Pres}}^{#1}_{#2}{#3}\setcounter{equation}{0}}
\def\ns{$n$-star}\setcounter{equation}{0}
\def\nt{$n$-tilting}\setcounter{equation}{0}
\def\Ht#1#2#3{{{\mathrm{Hom}}_{#1}({#2},{#3})}\setcounter{equation}{0}}
\def\qp#1{{${(#1)}$-quasi-projective}\setcounter{equation}{0}}
\def\mr#1{{{\mathrm{#1}}}\setcounter{equation}{0}}
\def\mc#1{{{\mathcal{#1}}}\setcounter{equation}{0}}
\def\ms#1{{{\mathscr{#1}}}\setcounter{equation}{0}}
%%%%%%%%FROM latexexam.tex
%\theoremstyle{definition}
\newcommand{\LL}{\ell\ell\,}
\newcommand{\edge}{\ar@{-}}
\newcommand{\bxy}{\xymatrix}
\newtheorem{Th}{Theorem}[section]%[section]%
\newtheorem{Def}[Th]{Definition}
\newtheorem{Lem}[Th]{Lemma}
\newtheorem{Pro}[Th]{Proposition}
\newtheorem{Cor}[Th]{Corollary}
\newtheorem{Rem}[Th]{Remark}
\newtheorem{Exm}[Th]{Example}
\newtheorem{Obs}[Th]{Observation}
\def\Pf#1{{\noindent\bf Proof}.\setcounter{equation}{0}}
\def\>#1{{ $\Rightarrow$ }\setcounter{equation}{0}}
\def\<>#1{{ $\Leftrightarrow$ }\setcounter{equation}{0}}
\def\bskip#1{{ \vskip 20pt }\setcounter{equation}{0}}
\def\sskip#1{{ \vskip 5pt }\setcounter{equation}{0}}
\def\zskip#1{{ \vskip 10pt }\setcounter{equation}{0}}
\def\bg#1{\begin{#1}\setcounter{equation}{0}}
\def\ed#1{\end{#1}\setcounter{equation}{0}}
\def\r#1{{\rm{#1}}}
\def\adr{\mr{add}_{\mc{D}^b(\mr{mod}R)}}
\def\ads{\mr{add}_{\mc{D}^b(\mr{mod}S)}}
\def\kbr{\ms{K}^b(\mc{P}_R)}
\def\kbri{\ms{K}^b(\mc{I}_R)}
\def\kbs{\ms{K}^b(\mc{P}_S)}
\def\dbr{\ms{D}^b(R)}
\def\dbs{\ms{D}^b(S)}
\def\dr#1{\ms{D}^{#1}(R)}
\def\ds#1{\ms{D}^{#1}(S)}
\def\kfr{\ms{K}^{-,b}(\mc{P}_R)}
\def\kfri{\ms{K}^{+,b}(\mc{I}_R)}
\def\kfs{\ms{K}^{-,b}(\mc{P}_S)}
\def\ks#1{\ms{K}^{#1}(\mc{P}_S)}
\def\kr#1{\ms{K}^{#1}(\mc{P}_R)}
\def\kri#1{\ms{K}^{#1}(\mc{I}_R)}
\def\Ok#1{\Omega_{\ms{D}}^{#1}}
\def\Oki#1{\Omega^{\ms{D}}_{#1}}
\def\Hr{\mr{Hom}_{\ms{D}(R)}}
\def\Hs{\mr{Hom}_{\ms{D}(S)}}

%\def\adr{\mr{add}_{\mc{D}^b(\mr{mod}R)}}
%\def\ads{\mr{add}_{\mc{D}^b(\mr{mod}S)}}

%%%%%%%%%%

%%%%%%%%%%%%%%%%%%%%%%%%%%%%%%%%%%%%%%%%%%%%%%%%%%%%%%%%%%%%%%%%%%%%%%%%%%%%%%%%%
%**************************±êÌâ¡¢ÕªÒª¡¢·ÖÀàºÅ¡¢¹Ø¼ü×Ö**************************

\title{\bf Derived categories and syzygies}
\smallskip
\author{{\small Jiaqun WEI} {\thanks {Supported by the National Science Foundation of China
(No.10971099)}}%\\
%\small Department of Mathematics,
%Nanjing Normal University \\
%\small Nanjing 210097, P.R.China\\ \small Email:
%weijiaqun@njnu.edu.cn}
}
\date{}
\maketitle
\baselineskip 17pt%16.5pt%18pt%17.2pt%20pt%19pt%14pt%\baselineskip  25.5pt
%%%%%%%\hskip
%
% Abstract ------------------------------------------------------
%
\begin{abstract}
\vskip 10pt%
We introduce syzygies for derived categories and study their
properties. Using these, we prove the derived invariance of the
following classes of artin algebras: (1) syzygy-finite algebras,
(2) Igusa-Todorov algebras, (3) AC algebras, (4) algebras
satisfying the finitistic Auslander conjecture, and (5) algebras
satisfying the generalized Auslander-Reiten conjecture. In
particular, Gorenstein CM-finite algebras are derived invariants.
\zskip\
%
%\noindent 2000 Mathematics Subject Classification:

\noindent MSC2000: \ \  Primary 18E30 16E05 Secondary 18G35 16G10
16E65
%
%\sskip\
%

\noindent {\it Keywords}: \ \ syzygy; derived category; derived
equivalence; Igusa-Todorov algebra; CM-finite algebra

\hskip 40pt \ Auslander condition
\end{abstract}
%\smallskip
%
\vskip 10pt
%
% ----------------------------------------------------------------------
%% ----------------------------------------------------------------------
%\def\baselinestretch{1}

%½éÉÜµ¼³ö·¶³ëµÄ×÷ÓÃ£¬Ó¦ÓÃ·¶Î§£¬¼òµ¥ÌáºÏ³åÔÚÍ¬µ÷´úÊýÖÐµÄ³£¼ûÐÔ£¬³£ÓÃÐÔ¡£
%Derived category.
%%%%%%%%%%%%%%%%%%%%%%%%%%%%%%%%%%%%%%%%%%%%%%%%%%%%%%%%%%%%%%%%%%%%
%
\vskip 30pt
\section{Introduction}

\hskip 15pt
%%%%%%%%%%%%%%%%%%%%%%%%%%%%%%%%%%%%%%%%%%%%%%%%%%%%%%%%%%%%%%%%%%%%

Syzygies were introduced by Hilbert [\ref{Hi}] in 1890. Nowadays,
syzygies are of general importance in algebraic geometry, homology
algebra and representation theory of groups and algebras etc.. The
advantage of syzygies is that they contain important information
of modules and they can test some homological dimensions. For
instance, Zimmermann-Huisgen and her coauthors
% provided the first counterexample to
%the first finitistic dimension conjecture by studying the
%structure of syzygies in monomial algebras [\ref{Zhi}] and it was
 showed that the structure of syzygies is important and powerful in
understanding and compute various finitistic dimension of some
algebras such as string algebras, monomial algebras and serial
algebras [\ref{GZ}, \ref{Zhi}, \ref{Zhs}, \ref{Zho}].

Monomial algebras and serial algebras are all syzygy-finite. By
definition, syzygy-finite algebras means that there is an integer
$s$ such that the class of all $n$-th syzygies, where $n>s$, is
representation finite, or equivalently, the number of
non-isomorphic indecomposable modules in the class is finite.
Syzygy-finite algebras have nice homological properties. They
satisfy a group of homological conjectures related to the
finitistic dimension conjecture, for instance, Auslander
conjecture (c.f. [\ref{Hpc}]) and Auslander-Reiten conjecture
[\ref{ARc}] etc.. Recently, algebras of finite Cohen-Macaulay
type, or CM-finite algebras, are more attractive, see for instance
[\ref{LZ}] and references therein.  Note that syzygy-finite
algebras are clearly CM-finite. We don't know if CM-finite
algebras are also syzygy-finite, but for an important class of
algebras, that is, Gorenstein algebras, they are the same.

In this paper, we extend the notion of syzygy to derived
categories and study their properties. As is well known, derived
categories are very important in modern study of algebraic
geometry and representation theory of groups and algebras, and
derived equivalences play an important role in the study. In
particular, there are remarkable conjectures concerning derived
equivalences. For instance, one famous conjecture in algebraic
geometry, first made by Bondal and Orlov, asserts that if $X_1$
and $X_2$ are birational smooth projective Calabi-Yau varieties of
dimension $n$, then there is an equivalence between their derived
categories [\ref{BO}]. While in representation theory of groups,
the famous Abelian Defect Group Conjecture of Brou\'{e} claims
that a block algebra $A$ of a finite group algebra and its Brauer
correspondent $B$ is derived equivalent provided that their common
defect group is abelian, see for instance [\ref{Rs}].% Derived
%equivalences preserve many significant properties, such as the
%number of non-isomorphic simple modules, the Hochschild cohomology
%groups, the finiteness of finitistic and global dimensions, and
%Gorenstein algebras etc..

Syzygies in derived category is powerful concerning derived
equivalences, as we show in this paper. For example, using
syzygies in derived category, we prove that derived equivalences
preserve the syzygy-finiteness of algebras. This provides an
important way to obtain syzygy-finite algebras, in particular,
CM-finite algebras. The reader is referred to [\ref{Ws}] for other
ways to obtain syzygy-finite algebras. Moreover, using syzygies in
derived category, we also prove that derived equivalences preserve
the following interesting classes of algebras: Igusa-Todorov
algebras, AC-algebras, algebras satisfying the finitistic
Auslander conjecture and algebras satisfying the
generalized Auslander-Reiten conjecture.% This gives an important
%property of the above mentioned algebras.

Igusa-Todorov algebras were introduced in [\ref{Wit}] in
connection with the study of the finitistic dimension conjecture
using Igusa-Todorov functor. Such algebras have finite finitistic
dimension. The class of Igusa-Todorov algebras is large, including
syzygy-finite algebras, algebras with representation dimension at
most three, algebras with radical cube zero and most algebras
which were recently proved to have finite finitistic dimension,
see [\ref{Wit}] for details. We refer to [\ref{Wit}, \ref{Ws}] for
other methods to judge when an
algebra is Igusa-Todorov.% an Iguan-Todorov algebra.

AC-algebras are algebras satisfying Auslander's condition. They
are studied in detail in [\ref{CH}] recently. Auslander's
condition (AC) for an algebra $R$ asserts that for every finitely
generated left $R$-module $M$ there is an integer $n = n_M$,
called Auslander bound of $M$, such that $\mr{Ext}^i_R(M,N) = 0$
for all $i>n$, whenever $N$ is a finitely generated left
$R$-module satisfying $\mr{Ext}^i_R(M,N) = 0$ for all but finitely
many $i$. Auslander conjectured all finite dimensional algebras
satisfy Auslander's condition (c.f.[\ref{Hpc}]). However, the
conjecture fails by counterexamples firstly given in [\ref{JS}]. A
revisited version of Auslander conjecture, named the finitistic
Auslander conjecture, asserts that the finitistic Auslander bound
of every algebra is finite [\ref{Wab}]. Note that the finitistic
Auslander conjecture implies the finitistic dimension conjecture.
In [\ref{Wab}], the generalized Auslander-Reiten conjecture is
also formulated, which asserts that for an algebra $R$, if $M$ is
a finitely generated left $R$-module such that
$\mr{Ext}_R^i(M,M\oplus R)=0$ for all $i\ge n$, then projective
dimension of $M$ is at most $n$. In the special case $n=1$, it is
just Auslander-Reiten conjecture [\ref{ARc}].

Though we work on artin algebras and finitely generated left
modules throughout this paper, syzygies (resp., cosyzygies)
defined here for derived categories make sense for the derived
category of any abelian categories with enough projective (resp.,
injective) objects. It is expected that these notions can have
nice applications in more areas.% in algebraic geometry and representation
%theory of groups and algebras.

The paper is organized as follows. In Section 2 we introduce
notations used in the paper and recall some basic facts on derived
categories. We introduce syzygies for derived categories and study
their basic properties in Section 3. In Sections 4 and 5, we prove
that derived equivalences preserve syzygy-finite algebras and
Igusa-Todorov algebras respectively. In Section 6, we prove that
AC-algebras, algebras satisfying the finitistic dimension
Auslander conjecture and algebras satisfying the generalized
Auslander-Reiten conjecture also have derived invariance.
% Some examples are presented in Section 6.

%%%%%%%%%%%%%%%%%%%%%%%%

\zskip\

\section{Preliminaries} %Derived categories and derived equivalences}

\hskip 17pt In this section, we recall some basic definitions and
facts which are necessary for our proofs.

% Now we introduce notations used in this paper.
Let $R$ be an artin
algebra, which means $R$ is a finitely generated $A$-module over a
commutative artin ring $A$. We denote by $\mr{mod}R$ the category
of finitely generated left $R$-modules. The subcategory of
$\mr{mod}R$ consisting of projective (resp., injective) modules is
denoted by $\mc{P}_R$ (resp., $\mc{I}_R$). We denote $fg$ the
composition of homomorphisms $f: L\to M$ and $g: M\to N$.
%
% All subcategories are full.

 We strengthen that we work on {\bf
chain} complexes. Let $\mc{C}$ be a class of $R$-modules.
 A (chain) complex $X$ over $\mc{C}$ is a set $\{X_i\in\mc{C} ,  i\in{\bf Z}\}$ equipped
with a set of homomorphisms $\{d_X^i: X_{i}\to X_{i-1}, i\in{\bf
Z} | d_X^{i+1}d_X^i=0\}$. We usually write $X=\{X_i,d_X^i\}$. A
chain map $f$ between complexes, say from $X=\{X_i,d_X^i\}$ to
$Y=\{Y_i,d_Y^i\}$, is a set of maps $f=\{f_i: X_i\to Y_i\}$ such
that $f_id_Y^i=d_X^if_{i-1}$.  A complex $X=\{X_i,d_X^i\}$ is
right (resp., left) bounded if $X_i=0$ for all but finitely many
negative (resp., positive) integers $i$. A complex $X$ is bounded
if it is both left and right bounded, equivalently, $X_i=0$ for
all but finitely many $i$. A complex $X$ is homologically bounded
if all but finitely many homologies of $X$ are zero.

Let $I\subseteq {\bf Z}$ be an %$X=\{X_i,d_X^i\}$ be a complex and
interval. We say that a complex $X=\{X_i,d_X^i\}$ has a
homological (resp., representation) interval $I$, if the
homologies $\mr{H}^i(X)=0$ (resp., terms $X_i=0$) for all
$i\not\in I$. We identify an $R$-module with a complex concentred
on the 0-th term, i.e., a complex with representation interval
$[0,0]$.

%$\mathscr{ABCDK} \mathfrak{ABCD,abcd}$

The category of all complexes over $\mc{C}$ with chain maps is
denoted by $\ms{C}(\mc{C})$. The homotopical category of complexes
over $\mc{C}$ is denoted by $\ms{K}(\mc{C})$. When $\mc{C}$ is an
abelian category, then the derived category of complexes over
$\mc{C}$ is well defined and is denoted by $\ms{D}(\mc{C})$. The
subcategories of $\ms{K}(\mc{C})$ and $\ms{D}(\mc{C})$ consisting
of bounded (resp., right bounded, left bounded) complexes are
denoted by $\ms{K}^b(\mc{C})$ (resp.,
 $\ms{K}^{-,b}(\mc{C})$, $\ms{K}^{+,b}(\mc{C})$) and
$\ms{D}^b(\mc{C})$ respectively. Let $I$ be an interval.
 We denote by $\ms{K}^{\mr{r}I}(\mc{C})$ (resp.,
 $\ms{K}^{\mr{h}I}(\mc{C})$)
and $\ms{D}^{I}(\mc{C})$ the subcategory of $\ms{K}(\mc{C})$ and
$\ms{D}(\mc{C})$ consisting of complexes with representation
(resp., homological) interval $I$, respectively. Similarly, we
denote by $\ms{D}^{I}(\mc{C})$ the subcategory of $\ms{D}(\mc{C})$
consisting of complexes with  homological interval $I$. For
instance, $\kr{\mr{r}[0,0]}$ is just the subcategory $\mc{P}_R$
and $\dr{[0,0]}$ is just the same as $\mr{mod}R$.

%We also denote by $\kr{\mr{-,h}I}$ (resp., $\kr{\mr{r}I}$) the
%subcategory of $\kfr$ consisting of complexes with homological
%(resp., representation) interval $I$. We also have similar notions
%$\kri{\mr{+,h}I}$ and $\kri{\mr{r}I}$.

We simply write  $\ms{D}(R)$ and  $\dbr$ for $\ms{D}(\mr{mod}R)$
and $\ms{D}^b(\mr{mod}R)$ respectively. It is well known that
$\dbr\subseteq\ms{D}(R)$, $\kbr\subseteq\kfr$,
$\ms{K}^{b}(\mc{I}_R)\subseteq\ms{K}^{+,b}(\mc{I}_R)$ are all
triangulated categories. For basic results in triangulated
categories and derived categories, we refer to [\ref{Hpt}] and
[\ref{V}].
 Moreover, $\dbr$ is equivalent to both
$\kfr$ and $\ms{K}^{+,b}(\mc{I}_R)$ as triangulated categories.

We denote by $[-]$ the shift functor on complexes. In fact, for a
complex $X$, the complex $X[1]$ is obtained from $X$ by shifting
$X$ to the left one degree. The notation $\mr{add}X$ denotes the
class of all direct summands of finite sums of a complex $X$.

Two algebras $R$ and $S$ are said to be derived-equivalent if
$\dbr$ and $\dbs$ are equivalent as triangulated categories.
Rickard [\ref{Ric}] proved that two algebras $R$ and $S$ are
derived-equivalent if and only if there is a complex $T\in\kbr$
with $S\simeq \Hr(T,T)$ such that $\Hr(T,T[i])=0$ for all $i\neq
0$ and $\mr{add}T$ generates $\kbr$ as a triangle category. Such
complex is called a tilting complex. In fact, if $\mc{F}:\ \dbr
\rightleftarrows\ \dbs\ :\mc{G}$ define an equivalence, then
$T:=\mc{G}(S)$ is just a tilting complex. In particular, if $T$ is
a tilting $R$-module, then $T$ induces a derived equivalence
between $\ms{D}^b(R)$ and $\ms{D}^b(\mr{End}_RT)$ [\ref{Hpt}]. We
refer to [\ref{HX}] for recent development on constructing derived
equivalences.

%Recall that, for a complex $M\in\dbr$, a projective resolution of
%$M$ is a complex $P\in\kfr$ such that $P\simeq M$. In case $M$ is
%an $R$-module, $P$ is just a usual projective resolution of $M$.
%
%We often consider these as identifications.

%Let $M,N\in\dbr$. The notion $\Hr(M,N)$ always is taken in

For an interval $I\subseteq {\bf Z}$, we denote by $\sigma_{I}(X)$
the brutal truncated complex which is obtained from the complex
$X$ by replacing each $X_i$, where $i\not\in I$, with 0. For
instance, if $M$ is a complex in $\kr{\mr{r}[n,\infty)}$, then
$\sigma_{[n,\infty)}(M)=M$.

Throughout the paper, $R$ stands for an artin algebra.
Homomorphisms and isomorphisms between complexes always means
in $\ms{D}(R)$. %the corresponding triangulated category.
 In
particular, restricting to module category, they are just the
usual homomorphisms and isomorphisms in module category.

%
%For basic results in triangulated categories and derived
%categories, we refer to [\ref{Hpt}] and [\ref{V}].

%%%%%%%%%%%%%%%%%%%%%%%%%%%%%%%%%%%%%%%%%%%%%%%%%%%%%%%%%%%%%%%%%%%%
%
\vskip 30pt
\section{Syzygies in derived categories}

\hskip 17pt %We first introduce the notion of syzygies in derived
%category as follows.

Recall that  $\dbr$ is equivalent to $\kfr$ as triangulated
categories. For a complex $M\in\dbr$, a projective resolution of
$M$ is a complex $P\in\kfr$ such that $P\simeq M$ in $\ms{D}(R)$.
In case $M$ is an $R$-module, $P$ is just a usual projective
resolution of $M$.

\bg{Def}\label{Sd}%Syzygy definition
Let $M\in\dbr$ and $n\in {\bf Z}$. Let $P$ be  a projective
resolution of $M$. We say that a complex in $\dbr$ is an $n$-th
syzygy of $M$, if it is isomorphic to
$(\sigma_{[n,\infty)}(P))[-n]$ in $\ms{D}(R)$. In the case, the
$n$-th syzygy of $M$ is denoted by $\Ok{n}(M_P)$, or simply by
$\Ok{n}(M)$ if there is no danger of confusion.
% if we don't care the
%projective resolution.
%(or  if there is no danger of confusion)
\ed{Def}

%%%%%%%%%%%%%%%%%%%%%%%%%%%%%%%%%%%%%%%%%%%%%%%
%\zskip\

Thus, the $n$-th syzygy of a complex $M\in \dbr$ depends on the
choice of projective resolution of  $M$. In this case that $M$ is
an $R$-module and $n>0$, the brutal truncated complex
$(\sigma_{[n,\infty)}(P))[-n]$ is just the projective resolution
of the $n$-th syzygy of $M$ defined by $P$. Hence, syzygies of $M$
defined here coincide with the usual syzygies in module category.

%We have the following easy observations.
%
%$(1)$  All syzygies are in $\kr{\mr{r}[0,\infty)}$, by definition.
%Hence, for any complex $M$ and integer $n$, we have that
%$\Hr(P,\Ok{n}(M)[i])=0$, whenever $P$ is projective and
%$i>0$.
%
%$(2)$ A direct sum of an $n$-th syzygy of $M$ and a projective
%module is still an $n$-th syzygy of $M$.
%
%
%$(3)$ $\Ok{n}(M[m])\simeq \Ok{n-m}(M)$ and
%$\Ok{n+m}(M)\simeq\Ok{m}(\Ok{n}(M))$, for any complex $M$ and
%integers $n,m$.
%
%$(4)$ If $M\in\dr{[a,b]}$ for some $a\le b$, then
%$\Ok{n}(M)\subseteq \dr{[0,b-n]}$ for any $n\in [a,b]$,
%$\Ok{n}(M)\subseteq \dr{[0,0]}=\mr{mod}R$ for any $n\in
%[b,\infty]$ and $\Ok{n}(M)\subseteq \dr{[0,b-a]}$ for any $n\in
%[a,\infty]$. In particular, $\Ok{b}(M)$ is an $R$-module and
%$M=\Ok{a}(M)[a]$.
%
%$(5)$ A complex in $\dbr$ is in $\kbr$ if and only if any/some
%syzygy of $M$ is also in $\kbr$.
%
%$(6)$  $\Oki{n}(M\oplus N)=\Oki{n}(M)\oplus \Oki{n}(N))$, for any
%complexes $M,N$ and any integer $n$.

We leave to the reader the entire proof of the following lemma.% since it is
%trivial.

\bg{Lem}\label{So}%Syzygy Obsversation
Let $M,N\in\dbr$ and $n,m,a,k,i$ be integers. Let  $P$ be a
projective resolutions of $M$.

$(1)$  $\Ok{n}(M)\in\dr{[0,\infty)}$ and $\Hr(Q,\Ok{n}(M)[i])=0$
for any projective module $Q$  and any integer $i>0$.

$(2)$ If $M\in\dr{[a,k]}$ for some $a\le k$, then $\Ok{n}(M)$ is
contained in $\dr{[0,0]}$, $\dr{[0,k-n]}$, $\dr{[a-n,k-n]}$ for
cases $n\ge k$, $a\le n\le k$, $n\le a$, respectively.
%
%\hskip 20pt (i) $\Ok{n}(M)\subseteq \dr{[0,b-n]}$, for any $n\in
%[a,b]$,
%
%\hskip 20pt (ii) $\Ok{n}(M)\subseteq \dr{[0,0]}=\mr{mod}R$, for
%any $n\in [b,\infty]$ and
%
%\hskip 20pt (iii) $\Ok{n}(M)\subseteq \dr{[a-n,b-n]}$ for any
%$n\in [-\infty,a]$.
%
 In particular, $\Ok{k}(M)$ is isomorphic to an $R$-module and
$M\simeq\Ok{a}(M)[a]$. % in $\ms{D}^b(R)$.

$(3)$ $\Ok{n+m}(M[m])\simeq \Ok{n}(M)$.

$(4)$ $\Ok{n+m}(M)\simeq\Ok{m}(\Ok{n}(M))$, for $m\ge 0$.
%$\Ok{n+m}(M)=(\Ok{n}(M))[-m]$, for $m< 0$.

$(5)$ $\Ok{n}(M)\oplus Q$ is also an $n$-th syzygy of $M$, where
$Q$ is a projective module.

$(6)$ $M\in\kbr$ if and only if any/some syzygy of $M$ is also in
$\kbr$.

$(7)$  $\Ok{n}(M\oplus N)\simeq\Ok{n}(M)\oplus \Ok{n}(N))$.
% denoted by $N=\Ok{n}(M)$,
\ed{Lem}
%
%\Pf. (1) Clearly $(\sigma_{[n,\infty)}(P))[-n]
%\in\kr{\mr{r}[0,\infty)}$. So that $\Ok{n}(M)\simeq
%(\sigma_{[n,\infty)}(P))[-n]\in\dr{[0,\infty)}$ by
%definition. Note that %the homologies
% $\Hr(R,(\Ok{n}(M))[i])
% \simeq\mr{H}^i(\Ok{n}(M))=0$ for all $i>0$, so the conclusion follows.
%
%(2) If $M\in\dr{[a,k]}$, then $P\in\kr{\mr{h}[a,k]}\bigcap
%\kr{\mr{r}[a,\infty)}$. Hence,
%$\sigma_{[n,\infty)}(P)\in\kr{\mr{h}I}\bigcap
%\kr{\mr{r}[n,\infty)}$, where $I$ is the interval $[n,n], [n,k],
%[a,k]$ for corresponding cases $n\ge k$, $a\le n\le k$, $n\le a$
%respectively. It follows that $\Ok{n}(M)$ is contained in
%$\dr{[0,0]}$, $\dr{[0,k-n]}$ and $\dr{[a-n,k-n]}$ for cases $n\ge
%k$, $a\le n\le b$ and $n\le a$ respectively. Note that
%$\dr{[0,0]}=\mr{mod}R$, so $\Ok{k}(M)$ is isomorphic to an
%$R$-module. It is also obvious that $\Ok{a}(M)[a]\simeq M$.
%
%(3) Note that $P[m]$ is a projective resolution of $M[m]$ and that
%$\sigma_{[n+m,\infty)}(P[m])=(\sigma_{[n,\infty)}(P))[m]$. Hence
%$\Ok{n+m}(M[m])\simeq \Ok{n}(M)$.
%
%(4) Note that
%$\sigma_{[n+m,\infty)}(P)=\sigma_{[n+m,\infty)}(\sigma_{[n,\infty)}(P))$
%for $m\ge 0$, so we have that $\Ok{n+m}(M) \\\simeq
%\Ok{n+m}((\Ok{n}(M))[n])\simeq \Ok{m}(\Ok{n}(M))$ by (3).
%
%
%(5) Let $P'=\{P'_i,d_{P'}^i\}$, where $P'_n=Q=P'_{n-1}$,
%$d_{P'}^n=1$ and $P'_i=0$ for other $i$. Then $P'\simeq 0$. Hence
%$P\oplus P'$ is also a projective resolution of $M$. It is easy to
%see that $\Ok{n}(M_P)\oplus Q\simeq\Ok{n}(M_{P\oplus P'})$.
%
%(6) and (7) are also easy.
%%
%\ \hfill$\Box$
%

%%%%%%%%%%%%%%%%%%%%%%%%%%%%%%%%%%%%%%%%%%%%%%%
\zskip\

%%%%%%%%%%%%%%%%%%%%%%%%%%%%%%%%%%%%%%%%%%%%%%%
\zskip\

Let $f$ be a chain map between complexes, say from
$X=\{X_i,d_X^i\}$ to $Y=\{Y_i,d_Y^i\}$. Recall that the cone of
$f$, denoted by $\mr{cone}(f)$, is a complex such that, for each
$i$, $(\mr{cone}(f))_i=Y_i\oplus X_{i-1}$, $d_{\mr{cone}(f)}^i=(\
{^{d_Y^i\ f_{i-1}}_{0\ \ \ d_X^{i-1}}})$. It is well known that
there is a canonical triangle $X\to^f Y\to \mr{cone}(f)\to$ in
$\ms{D}(R)$. Using the construction of cones, one obtains the
following canonical triangles in $\ms{D}(R)$,  for any complex $X$
and any $n\ge m$:

\bg{verse}  ($\sigma$1) $(\sigma_{[n+1,\infty)}(X))[-1]\to
\sigma_{[m,n]}(X)\to \sigma_{[m,\infty)}(X)\to$,

($\sigma$2) $(\sigma_{[m,n]}(X))[-1]\to
\sigma_{[-\infty,m-1]}(X)\to \sigma_{(-\infty,n]}(X)\to$, and

($\sigma$3) $(\sigma_{[n,\infty)}(X))[-1]\to X_{(-\infty,n-1]}\to
X\to$.\ed{verse}

%%%%%%%%%%%%%%%%%%%%%%%%%%%%%%%%%%%%%%%%%%%%%%%%%%%%%%%%%%%%%%%%%%%%
\zskip\

\bg{Lem}\label{Stri}%syztri
Let $M\in \dbr$. Then
% and $P$ be a projective resolution of $M$. Then
%there is a triangle $(\Ok{n+1}(M_P))[n-m]\to
%(\sigma_{[m,n]}(P))[-m]\to \Ok{m}(M_P)$, for any $n\ge m$. In
%particular, for any $n$, there is a triangle $\Ok{n+1}(M_P)\to
%P_n\to \Ok{n}(M_P)$, where $P_n$ is the $n$-th term of $P$.
%
there is a triangle $(\Ok{n+1}(M))[n-m]\to B\to \Ok{m}(M)\to$,
where $B\in\kr{\mr{r}[0,n-m]}$ and $n\ge m$. In particular, for
any $n$, there is a triangle $\Ok{n+1}(M)\to Q\to \Ok{n}(M)\to$
with $Q$ projective.
\ed{Lem}

\Pf. Let $P$ be a projective resolution of $M$. Then the first
triangle is obtained from $(\sigma 1)$ by shifting, %$m$ times,
letting $B=(\sigma_{[m,n]}(P))[-m]\in\kr{\mr{r}[0,n-m]}$) in the
case. Taking $m=n$, we then obtain the second triangle.
\ \hfill$\Box$

%%%%%%%%%%%%%%%%%%%%%%%%%%%%%%%%%%%%%%%%%%%%%%%
\zskip\

Recall that two modules $N,N'$ are projectively equivalent if
$N\oplus P\simeq N'\oplus P'$ for some projective modules $P,P'$.
Let $M$ be a module. It is well known that any two $n$-th syzygies
of $M$ are projectively equivalent. It is also the case for
syzygies in derived category. We say that two complexes $N,N'$ are
projectively equivalent if $N\oplus P\simeq N'\oplus P'$ for some
projective modules $P,P'$.

\bg{Pro}\label{Sp}%syzygy proj equiv
Let $M\in\dbr$ and $P,P'$ be two projective resolutions of $M$.
Then $\Ok{n}(M_P)$ and $\Ok{n}(M_{P'})$ are projectively
equivalent. %%%ÏÂÃæconverse¿ÉÄÜ²»¶Ô
%Conversely, if $N\in\dbr$ is projectively equivalent
%to $\Ok{n}(M_P)$, then $N$ is also an $n$-syzygy of $M$.
%
\ed{Pro}

\Pf. We prove by induction on $n$.
% Without loss of generality, we
%may assume that $P,P'\in\kr{\mr{r}[a,\infty)}$ for some integer
%$a$. Since $\sigma_{[n,\infty)}(P)\simeq M\simeq
%\sigma_{[n,\infty)}(P')$ for all $n\le a$, we obtain that
%$\Ok{n}(M_P)\simeq (\sigma_{[n,\infty)}(P))[-n]\simeq
%(\sigma_{[n,\infty)}(P'))[-n]\simeq\Ok{n}(M_{P'})$ in case $n\le
%a$.
Assume that $M\in\dr{[a,k]}$ for some integers $a<k$. Then
$\Ok{n}(M_P)\simeq M[-n]\simeq \Ok{n}(M_{P'})$ in case $n\le a$,
by Lemma \ref{So} (2).

Now we consider the case $n+1$.  By the induction assumption,
there are projective modules $Q,Q'$ such that $\Ok{n}(M_{P})\oplus
Q\simeq \Ok{n}(M_{P'})\oplus Q'$. Then we have triangles

\sskip\

\hskip 100pt$\Ok{n+1}(M_P)\to P_n\oplus Q\to \Ok{n}(M_P)\oplus
Q$,\hskip 20pt and

\hskip 100pt$\Ok{n+1}(M_{P'})\to {P'_n}\oplus Q'\to
\Ok{n}(M_{P'})\oplus Q'$,

\sskip\

\noindent where $P_n$ and ${P'_n}$ are projective, by Lemma
\ref{Stri}. Since we have a homomorphism between the last terms of
these two triangles and $\Hr(P_n\oplus Q,\Ok{n+1}(M_{P'})[1])=0$
by Lemma \ref{So} (1), there is the following commutative diagram:

\sskip\

 \setlength{\unitlength}{0.08in}
 \begin{picture}(50,12)
%
% \put(18,1){\makebox(0,0)[c]{$0$}}
% \put(27,1){\makebox(0,0)[c]{$0$}}
% \put(35,1){\makebox(0,0)[c]{$0$}}
%
%                 \put(18,3.4){\vector(0,-1){2}}
%                 \put(27,3.4){\vector(0,-1){2}}
%                 \put(35,3.4){\vector(0,-1){2}}

% \put(12,5){\makebox(0,0)[c]{$0$}}
%                             \put(14,5){\vector(1,0){2}}
 \put(18,3){\makebox(0,0)[c]{$\Ok{n+1}(M_{P'})$}}
                             \put(23,3){\vector(1,0){2}}
%                             \put(20,5){\vector(1,0){2}}
%                             \put(21,6){\makebox(0,0)[c]{$_f$}}
 \put(29,3){\makebox(0,0)[c]{${P'_n}\oplus Q'$}}
%                             \put(31,5){\vector(1,0){2}}
                             \put(34,3){\vector(1,0){2}}
%                             \put(32,6){\makebox(0,0)[c]{$_g$}}
 \put(43,3){\makebox(0,0)[c]{$\Ok{n}(M_{P'})\oplus Q'$}}
%                             \put(37,5){\vector(1,0){2}}
                             \put(49,3){\vector(1,0){2}}
% \put(41,5){\makebox(0,0)[c]{$0$}}

%                 \put(17.9,9){\line(0,-1){1}}
%                 \put(18.1,9){\line(0,-1){1}}
                 \put(18,7){\vector(0,-1){2}}
%                       \put(17,8){\makebox(0,0)[c]{$_{u_0}$}}
                 \put(29,7){\vector(0,-1){2}}
%                       \put(25,8){\makebox(0,0)[c]{$(_{{\ \theta\ }}^{{u_0}f})$}}
                 \put(43,7){\vector(0,-1){2}}
%                       \put(37,8){\makebox(0,0)[c]{$_{w_0}$}}

% \put(12,11){\makebox(0,0)[c]{$0$}}
%                             \put(14,11){\vector(1,0){2}}
 \put(18,9){\makebox(0,0)[c]{$\Ok{n+1}(M_P)$}}
                             \put(23,9){\vector(1,0){2}}
%                             \put(21,12){\makebox(0,0)[c]{${_{(1, 0)}}$}}
 \put(29,9){\makebox(0,0)[c]{$P_n\oplus Q$}}
                             \put(34,9){\vector(1,0){2}}
%                             \put(32,12.5){\makebox(0,0)[c]{$({_{1}^{0}})$}}
 \put(43,9){\makebox(0,0)[c]{$\Ok{n}(M_P)\oplus Q$}}
                             \put(49,9){\vector(1,0){2}}
%                             \put(33,10){\vector(-1,-1){4}}
%                             \put(30,8){\makebox(0,0)[c]{$_\theta$}}
% \put(41,11){\makebox(0,0)[c]{$0$}}
%
%                 \put(18,14.5){\vector(0,-1){2}}
%                 \put(27,14.5){\vector(0,-1){2}}
%                 \put(35,14.5){\vector(0,-1){2}}
%
%
%% \put(12,16){\makebox(0,0)[c]{$0$}}
%%                             \put(14,16){\vector(1,0){2}}
% \put(18,16){\makebox(0,0)[c]{$\Ok{n+1}(L)$}}
%%                             \put(20,16){\vector(1,0){2}}
%                             \put(22,16){\vector(1,0){2}}
% \put(27,16){\makebox(0,0)[c]{$M'$}}
%%                             \put(31,16){\vector(1,0){2}}
%                             \put(29,16){\vector(1,0){2}}
% \put(35,16){\makebox(0,0)[c]{$\Ok{n+1}(N)$}}
%%                             \put(37,16){\vector(1,0){2}}
%                             \put(39,16){\vector(1,0){2}}
% \put(41,16){\makebox(0,0)[c]{$0$}}

%                 \put(18,19){\vector(0,-1){2}}
%                 \put(27,19){\vector(0,-1){2}}
%                 \put(35,19){\vector(0,-1){2}}
%
%
% \put(18,21){\makebox(0,0)[c]{$0$}}
% \put(27,21){\makebox(0,0)[c]{$0$}}
% \put(35,21){\makebox(0,0)[c]{$0$}}

\end{picture}

Since the homomorphism in the right column is an isomorphism, we
have a canonical triangle

\sskip\

\hskip 100pt$\Ok{n+1}(M_P)\to P_n\oplus Q\oplus
\Ok{n+1}(M_{P'})\to {P'_n}\oplus Q'\to$.

\sskip\

\noindent  Note again that $P'_n\oplus Q'$ are projective and that
$\Hr(P'_n\oplus Q',\Ok{n+1}(M_P)[1])=0$ by Lemma \ref{So} (1), so
the above triangle splits. Hence we obtain that

\sskip\

\hskip 100pt$\Ok{n+1}(M_P)\oplus ({P'_n}\oplus Q')\simeq
(P_n\oplus Q)\oplus \Ok{n+1}(M_{P'})$.

\sskip\

\noindent It follows that $\Ok{n+1}(M_P)$ and $\Ok{n+1}(M_{P'})$
are projectively equivalent.
\ \hfill$\Box$

%%%%%%%%%%%%%%%%%%%%%%%%%%%%%%%%%%%%%%%%%%%%%%%
\zskip\

By the above result, we see that the $n$-th syzygy of $M$ is
unique up to projectively equivalences. By abuse of language, we
speak of the $n$-syzygy of a complex in $\dbr$.

We also have the following easy result. The proof is left to the
reader.

%%%%%%%%%%%%%%%%%%%%%%%%%%%%%%%%%%%%%%%%%%%%%%%%%%%%%%%%%%%%%%%%%%%%

\bg{Pro}\label{St}%syztri
{\rm (Dimension shifting)} Let $M\in \dbr$.
%
%$(1)$ There is a triangle $(\Ok{n+1}(M))[n-m]\to B\to \Ok{m}(M)$
%for some $B\in\kr{\mr{r}[0,n-m]}$, where $n\ge m$. In particular,
%for any $n$, there is a triangle $\Ok{n+1}(M)\to P\to \Ok{n}(M)$
%with $P$ projective.
%
%$(2)$
Assume that $N\in\dr{[c,d]}$ for some integers $c\le d$. Then
there is an isomorphism

\sskip\

\hskip 100pt$\Hr(\Ok{n+m}(M),N[j])\simeq \Hr(\Ok{n}(M),N[j+m])$

\sskip\

\noindent
 for
any integers $n,m,j$ such that $m\ge 1$ and $j>-c$.
\ed{Pro}
%
%\Pf. % (1) follows Lemma \ref{Stri}.% and Proposition \ref{Sp}. In
%%%fact, the desired triangle is obtained from the one in Lemma
%%%\ref{Stri} by correcting with some projective module or its shift.
%%
%%(2)
%Note that $N[-c]\in \dr{[0,d-c]}$, so $\Hr(P,N[j])=0$ for all
%projective $P$ and all $j>-c$. Hence, the isomorphism follows, by
%applying the functor $\Hr(-,N)$ to triangles $\Ok{n+i+1}(M)\to
%P_i\to \Ok{n+i}(M)\to$ obtained from Lemma \ref{Stri} for $0\le
%i<m$, where each $P_i$ is projective.
%%
%\ \hfill$\Box$

%
%%%%%%%%%%%%%%%%%%%%%%%%%%%%%%%%%%%%%%%%%%%%%%%%
\zskip\

The following result provides a way to compare syzygies for
complexes in some triangles.

\bg{Lem}\label{syzl}%
Let $M\to B\to N\to $ be a triangle in $\dbr$. If
$B\in\kr{\mr{r}(-\infty,k]}$ for some integer $k$, then
$\Ok{k}(M)\simeq\Ok{k+1}(N)$. Hence, $\Ok{n}(M)\simeq\Ok{n+1}(N)$
for all $n\ge k$.
\ed{Lem}

\Pf. Let $i: M\to B$ be the homomorphism in the triangle. Assume
that $P$ is a projective resolution of $M$, then we have a
homomorphism $f: P\to B$ in $\kfr$ followed from the equivalence
between $\dbr$ and $\kfr$.
% We can consider $M,B$ in $\kfr$ and
%then $f$ is a homomorphism between complexes.
 Note that $\mr{cone}(f)\in\kfr$ and $N\simeq
\mr{cone}(f)$, i.e., $\mr{cone}(f)$ is a projective resolution of
$N$.
%,
% so it is equivalent to show that
%$\Ok{b}(M)\simeq\Ok{b+1}(\mr{cone}(f))$.
 Since
$B\in\kr{\mr{r}(-\infty,k]}$, it is easy to see that
$\sigma_{[k+1,\infty)}(\mr{cone}(f))\simeq
\sigma_{[k+1,\infty)}(P[1])$, by the construction of
$\mr{cone}(f)$. It follows that
\sskip\
%
%\hskip 70pt$\Ok{b+1}(\mr{cone}(f))\simeq
\hskip 98pt$\Ok{k+1}(N)\simeq
(\sigma_{[k+1,\infty)}(\mr{cone}(f)))[-(k+1)]$

\hskip 145pt $\simeq(\sigma_{[k+1,\infty)}(P[1]))[-(k+1)]$

\hskip 145pt $\simeq\Ok{k+1}(M[1])\simeq\Ok{k}(M)$.
\sskip\

\noindent by Lemma \ref{So} (3). Hence, the conclusion follows.
%$\Ok{a-1}(M)\simeq\Ok{a-1}(\mr{cone}(f))$.
%
\ \hfill$\Box$

%
%%%%%%%%%%%%%%%%%%%%%%%%%%%%%%%%%%%%%%%%%%%%%%%%
\zskip\

For general triangles, we have the following result.

\bg{Pro}\label{Stn}%
Let $L\to M\to N\to $ be a triangle in $\dbr$. Then, for any
integer $n$, there exists a triangle
$\Ok{n}(L)\to\Ok{n}(M)\to\Ok{n}(N)\to$.
\ed{Pro}

\Pf. The proof is given by induction on $n$. We may assume that
$L,M,N\in\dr{[a,k]}$ for some integers $a\le k$. In case $n\le a$,
we obtain that $\Ok{n}(L)\simeq L[-n]$, $\Ok{n}(M)\simeq M[-n]$
and $\Ok{n}(N)\simeq N[-n]$, by Lemma \ref{So} (2). Hence we have
a triangle $\Ok{n}(L)\to\Ok{n}(M)\to\Ok{n}(N)\to$ by assumption.

Now consider the case $n+1$. By Lemma \ref{Stri}, we have
triangles $\Ok{n+1}(L)\to B\to \Ok{n}(L)\to$ and $\Ok{n+1}(N)\to
C\to \Ok{n}(N)\to$, where $B,C$ are projective. By the induction
assumption, there is a triangle
$\Ok{n}(L)\to\Ok{n}(M)\to\Ok{n}(N)\to$. Note that
$\Hr(C,(\Ok{n}(L))[1])=0$ by Lemma \ref{So} (1), so we can obtain
the following triangle commutative diagram for some $M'\in\dbr$.

\sskip\

 \setlength{\unitlength}{0.09in}
 \begin{picture}(50,19)
%
% \put(18,1){\makebox(0,0)[c]{$0$}}
% \put(27,1){\makebox(0,0)[c]{$0$}}
% \put(35,1){\makebox(0,0)[c]{$0$}}
%
                 \put(18,3.4){\vector(0,-1){2}}
                 \put(27,3.4){\vector(0,-1){2}}
                 \put(35,3.4){\vector(0,-1){2}}

% \put(12,5){\makebox(0,0)[c]{$0$}}
%                             \put(14,5){\vector(1,0){2}}
 \put(18,5){\makebox(0,0)[c]{$\Ok{n}(L)$}}
                             \put(21,5){\vector(1,0){2}}
%                             \put(20,5){\vector(1,0){2}}
%                             \put(21,6){\makebox(0,0)[c]{$_f$}}
 \put(27,5){\makebox(0,0)[c]{$\Ok{n}(M)$}}
%                             \put(31,5){\vector(1,0){2}}
                             \put(30,5){\vector(1,0){2}}
%                             \put(32,6){\makebox(0,0)[c]{$_g$}}
 \put(35,5){\makebox(0,0)[c]{$\Ok{n}(N)$}}
%                             \put(37,5){\vector(1,0){2}}
                             \put(39,5){\vector(1,0){2}}
% \put(41,5){\makebox(0,0)[c]{$0$}}

%                 \put(17.9,9){\line(0,-1){1}}
%                 \put(18.1,9){\line(0,-1){1}}
                 \put(18,9){\vector(0,-1){2}}
%                       \put(17,8){\makebox(0,0)[c]{$_{u_0}$}}
                 \put(27,9){\vector(0,-1){2}}
%                       \put(25,8){\makebox(0,0)[c]{$(_{{\ \theta\ }}^{{u_0}f})$}}
                 \put(35,9){\vector(0,-1){2}}
%                       \put(37,8){\makebox(0,0)[c]{$_{w_0}$}}

% \put(12,11){\makebox(0,0)[c]{$0$}}
%                             \put(14,11){\vector(1,0){2}}
 \put(18,11){\makebox(0,0)[c]{$B$}}
                             \put(20,11){\vector(1,0){2}}
%                             \put(21,12){\makebox(0,0)[c]{${_{(1, 0)}}$}}
 \put(27,11){\makebox(0,0)[c]{$B\oplus C$}}
                             \put(31,11){\vector(1,0){2}}
%                             \put(32,12.5){\makebox(0,0)[c]{$({_{1}^{0}})$}}
 \put(35,11){\makebox(0,0)[c]{$C$}}
                             \put(37,11){\vector(1,0){2}}
%                             \put(33,10){\vector(-1,-1){4}}
%                             \put(30,8){\makebox(0,0)[c]{$_\theta$}}
% \put(41,11){\makebox(0,0)[c]{$0$}}

                 \put(18,14.5){\vector(0,-1){2}}
                 \put(27,14.5){\vector(0,-1){2}}
                 \put(35,14.5){\vector(0,-1){2}}

% \put(12,16){\makebox(0,0)[c]{$0$}}
%                             \put(14,16){\vector(1,0){2}}
 \put(18,16){\makebox(0,0)[c]{$\Ok{n+1}(L)$}}
%                             \put(20,16){\vector(1,0){2}}
                             \put(22,16){\vector(1,0){2}}
 \put(27,16){\makebox(0,0)[c]{$M'$}}
%                             \put(31,16){\vector(1,0){2}}
                             \put(29,16){\vector(1,0){2}}
 \put(35,16){\makebox(0,0)[c]{$\Ok{n+1}(N)$}}
%                             \put(37,16){\vector(1,0){2}}
                             \put(39,16){\vector(1,0){2}}
% \put(41,16){\makebox(0,0)[c]{$0$}}

%                 \put(18,19){\vector(0,-1){2}}
%                 \put(27,19){\vector(0,-1){2}}
%                 \put(35,19){\vector(0,-1){2}}
%
%
% \put(18,21){\makebox(0,0)[c]{$0$}}
% \put(27,21){\makebox(0,0)[c]{$0$}}
% \put(35,21){\makebox(0,0)[c]{$0$}}

\end{picture}

Consider the triangle $M'\to B\oplus C\to \Ok{n}(M)\to$ from the
middle column in the diagram. Since $B\oplus C$ is a projective
$R$-module, $B\oplus C\in\kr{\mr{r}[0,0]}$. By Lemmas \ref{syzl}
and \ref{So}, we obtain that $\Ok{0}(M')\simeq
\Ok{1}(\Ok{n}(M))\simeq\Ok{n+1}(M)$. Note that $M'\in
\dr{[0,\infty)}$ followed from the triangle $\Ok{n+1}(L)\to M'\to
 \Ok{n+1}(N)\to$ from the up row in the diagram, so $\Ok{0}(M')\simeq M'$ by Lemma \ref{So} (2). It follows that
 $M'\simeq\Ok{n+1}(M)$ and hence we have a triangle $\Ok{n+1}(L)\to \Ok{n+1}(M)\to
 \Ok{n+1}(N)\to$.
\ \hfill$\Box$

%%%%%%%%%%%%%%%%%%%%%%%%%%%%%%%%%%%%%%%%%%%%%%%%%%%%%%%%%%%%%%%%%%%%
\zskip\

There is also the notion dual to syzygies. Recall that $\dbr$ is
equivalent to $\kfri$ as triangulated categories. For a complex
$M\in\dbr$, an injective resolution of $M$ is a complex
$I\in\kfri$ such that $I\simeq M$. In case $M$ is an $R$-module,
$I$ is just a usual injective resolution of $M$. Using the
injective resolution of a complex in $\dbr$, we can define the
notion of cosyzygies.

\zskip\

\noindent {\bf Definition \ref{Sd}'}\ \ %
{\it  Let $M\in\dbr$ and $n\in {\bf Z}$. Let $I$ be an injective
resolution of $M$. We say that a complex in $\dbr$ is an $n$-th
cosyzygy of $M$, if it is isomorphic to
$(\sigma_{(-\infty,n]}(I))[-n]$ in $\ms{D}(R)$. In the case, the
$n$-th cosyzygy of $M$ is denoted by $\Oki{n}(M_I)$, or simply by
$\Oki{n}(M)$ if there is no danger of confusion.}
% denoted by $N=\Ok{n}(M)$,
%
\zskip\
We state the dual results without proofs.
\zskip\
\noindent {\bf Lemma \ref{So}'}\ \  %Syzygy Obsversation
Let $M,N\in\dbr$ and $n,m,a,b,i$ be integers.% Let  $P$ be a
%injective resolutions of $M$.

$(1)$  $\Oki{n}(M)\in\dr{(-\infty,0]}$ and
$\Hr(\Oki{n}(M),Q[i])=0$ for any injective module $Q$  and any
integer $i>0$.

$(2)$ If $M\in\dr{[a,k]}$, for some $a\le k$, then $\Oki{n}(M)$ is
contained in $\dr{[0,0]}$, $\dr{[a-n,0]}$, $\dr{[a-n,k-n]}$ for
corresponding cases $n\le a$, $a\le n\le k$, $n\ge b$
respectively. In particular, $\Oki{a}(M)$ is isomorphic to an
$R$-module and $M\simeq\Oki{k}(M)[k]$.

$(3)$ $\Oki{n+m}(M[m])\simeq \Oki{n}(M)$.

$(4)$ $\Oki{n+m}(M)\simeq\Oki{m}(\Oki{n}(M))$ for $m\le 0$.

$(5)$ $\Oki{n}(M)\oplus Q$ is also an $n$-th cosyzygy of $M$,
where $Q$ is any injective module.

$(6)$ $M\in\kbri$ if and only if any/some cosyzygy of $M$ is also
in $\kbri$.

$(7)$  $\Oki{n}(M\oplus N)\simeq\Oki{n}(M)\oplus \Oki{n}(N))$.
% denoted by $N=\Ok{n}(M)$,
%
\zskip\

\noindent {\bf Lemma \ref{Stri}'}\ \ \ \ %syztri
{\it Let $M\in \dbr$. Then there is a triangle $\Oki{n}(M)\to B\to
(\Oki{m-1}(M))[m-n]\to$ for some $B\in\kri{\mr{r}[m-n,0]}$, where
$n\ge m$. In particular, for any $n$, there is a triangle
$\Oki{n}(M)\to I\to \Oki{n-1}(M)\to$ with $I$ injective.}
\zskip\

\noindent {\bf Proposition \ref{Sp}'}\ \  %syzygy proj equiv
{\it Let $M\in\dbr$ and $I,I'$ be two injective resolutions of
$M$. Then $\Oki{n}(M_I)$ and $\Oki{n}(M_{I'})$ are injectively
equivalent.}
\zskip\
\noindent {\bf Proposition \ref{St}'}\ \  {\rm (Dimension shifting)} %syztri
{\it Assume that $N\in\dr{[c,d]}$, where $c\le d$. Then there is
an isomorphism $\Hr(N,(\Oki{n}(M))[j])\simeq
\Hr(N,(\Oki{n+m}(M))[j+m])$ for any integers $n,m,j$ such that
$m\ge 1$ and $j>-d$.}
\zskip\
\noindent {\bf Lemma \ref{syzl}'}\ \  %
{\it Let $M\to B\to N\to $ be a triangle in $\dbr$. If
$B\in\kri{\mr{r}[a,\infty)}$ for some integer $a$, then
$\Oki{a}(M)\simeq\Oki{a+1}(N)$. Hence,
$\Oki{n}(M)\simeq\Oki{n+1}(N)$ for all $n\le a$.}
\zskip\
\noindent {\bf Proposition \ref{Stn}'}\ \ %
 Let $L\to M\to N\to $ be a
triangle in $\dbr$. Then there is a triangle
$\Oki{n}(L)\to\Oki{n}(M)\to\Oki{n}(N)\to$, for any integer $n$.

\zskip\

Let us remark that one can define syzygies (resp., cosyzygies) in
the derived category of any abelian category with enough
projective (resp., injective) objects.
%
%\ed{Rem}

%%%%%%%%%%%%%%%%%%%%%%%%%%%%%%%%%%%%%%%%%%%%%%%%%%%%%%%%%%%%%%%%%%%%
\vskip 15pt

\bskip\

\section{Syzygy-finite algebras}%Derived invariance}%

%%%%%%%%%%%%%%%%%%%%%%%%%%%%%%%%%%%%%%%%%%%%%%%%%%%%%%%%%%%%%%%%%%%%
\hskip 15pt\hskip 15pt

Let $\mc{C}\subseteq \dbr$. For an integer $n$, we denote by
$\Ok{n}(\mc{C})$ the class of all $n$-th syzygies of complexes in
$\mc{C}$. The class $\mc{C}$ is representation-finite provided
that $\mc{C}\subseteq \mr{add}M$ for some $M\in\dbr$, or
equivalently, the number of non-isomorphic indecomposable direct
summands of objects in $\mc{C}$ is finite. It is easy to see that
if $\mc{E}\subseteq\mc{C}$ and $\mc{C}$ is representation-finite
then $\mc{E}$ is also representation-finite. We say that $\mc{C}$
is $n$-{\it syzygy-finite} provided that $\bigcup_{i\ge
n}\Ok{i}(\mc{C})$ is representation-finite. By that $\mc{C}$ is
syzygy-finite, we mean that $\mc{C}$ is $n$-syzygy-finite for some
$n$. It follows from Lemma \ref{So} (3) that $\mc{C}$ is
syzygy-finite if and only if $\mc{C}[m]$ is syzygy-finite for
any/some integer $m$. %In case that $\mc{C}$ has the property
%$\Ok{1}(\mc{C})\subseteq \mc{C}$, it is easy to see that $\mc{C}$
%is $n$-syzygy-finite if and only if $\Ok{n}(\mc{C})$ is
%representation-finite.
%Note that if $\mc{C}$ is
%representation-finite, then $\Ok{n}(\mc{C})$ is
%representation-finite for any $n$.

An artin algebra $R$ is called syzygy-finite if the class
$\mr{mod}R$ is syzygy-finite. The following algebras are known to
be syzygy-finite.

\zskip\

$\bullet$ Algebras of finite representation type
% (in fact $0$-syzygy-finite).

$\bullet$ Algebras of finite global dimension.

$\bullet$ Monomial algebras % (indeed $2$-syzygy-finite)
[\ref{Zhi}].

$\bullet$ Left serial algebras % (indeed $2$-syzygy-finite)
[\ref{Zhs}].

$\bullet$ Torsionless-finite algebras, %  (indeed $1$-syzygy-finite)
c.f. [\ref{Rg}], including:

\hskip 30pt $\circ$ Algebras $R$ with $\mr{rad}^nR = 0$ and
$A/\mr{rad}^{n-1}A$  representation-finite.

\hskip 30pt $\circ$ Algebras $R$ with $\mr{rad}^2R = 0$.

\hskip 30pt $\circ$ Minimal representation-infinite algebras.

\hskip 30pt $\circ$ Algebras stably equivalent to hereditary
algebras.

\hskip 30pt $\circ$ Right glued algebras and left glued algebras.

\hskip 30pt $\circ$ Algebras of the form $R/\mr{soc}R$ with $R$ a
local algebra of quaternion type.

\hskip 30pt $\circ$  Special biserial algebras.

$\bullet$ Algebras possessing a left idealized extension which is
torsionless-finite (indeed $2$-syzygy-finite) [\ref{Wit}].

$\bullet$ Algebras possessing an ideal $I$ of finite projective
dimension such that $I\mr{rad}R=0$ and $R/I$ is syzygy-finite
[\ref{Ws}].

\zskip\

The following result gives a characterization of syzygy-finite
algebras in term of derived category.

\bg{Th}\label{Sfd}%Syzygy-finite derived
An algebra $R$ is syzygy-finite if and only if $\dr{[a,k]}$ is
syzygy-finite for any/some integers $a\le k$.
\ed{Th}%

\Pf.  Assume that $\dr{[a,k]}$ is syzygy-finite for some integers
$a\le k$. Since $(\mr{mod}R)[a]\subseteq \dr{[a,k]}$,
$(\mr{mod}R)[a]$ is syzygy-finite and hence $\mr{mod}R$ is
syzygy-finite, i.e., $R$ is a syzygy-finite algebra.

On the other hand, assume that $R$ is syzygy-finite. For any
integers $a\le k$, we have $\Ok{b}(\dr{[a,k]})\subseteq \mr{mod}R$
by Lemma \ref{So} (2). Since $\mr{mod}R$ is syzygy-finite, we see
that $\Ok{k}(\dr{[a,k]})$ is syzygy-finite. It follows that
$\dr{[a,k]}$ is syzygy-finite from Lemma \ref{So} (4).
\ \hfill$\Box$

%%%%%%%%%%%%%%%%%%%%%%%%%%%%%%%%%%%%%%%%%%%%%%%%%%%%%%%%%%%%%%%%%%%%
\zskip\

\bg{Pro}\label{syzbj}%functor Lemma
Let $\mc{C}, \mc{E}\subseteq \dbr$ and $\mc{B}\subseteq
\kr{\mr{r}(-\infty,k]}$ for some $b$. Assume that, for any
$E\in\mc{E}$, there is a triangle $C_E\to B_E\to E\to $ in $\dbr$
with $C_E\in\mc{C}$ and $B_E\in\mc{B}$. If $\mc{C}$ is
syzygy-finite, then $\mc{E}$ is also syzygy-finite.
%
%$(2)$ Assume moreover $\mc{D}\subseteq \kr{[a,b]}$ for some fixed
%integers $a,b$. If $\mc{C}$ is syzygy-finite from some $c$, then
%$\mc{E}$ is syzygy-finite from some $e$, where $e$ depending only
%on $a$ and $c$.
%
\ed{Pro}%

\Pf. Assume that $\mc{C}$ is $m$-syzygy-finite, for some $m$. For
any $E\in\mc{E}$, consider the triangle $C_E\to B_E\to E\to$ in
$\dbr$ with $C_E\in\mc{C}$ and $B_E\in\mc{B}$. Since $B_E\in
\kr{\mr{r}(-\infty,k]}$ by assumptions, we have that
$\Ok{i}(C_E)\simeq\Ok{i+1}(E)$ for all $i\ge k$, by Lemma
\ref{syzl}. Hence, $\bigcup_{i\ge n+1}\Ok{i}(\mc{E})$ is
representation-finite, for some $n\ge k$, if and only if the class
 $\{\Ok{i}(C_E)\ |\ E\in\mc{E},
i\ge n\}$ is representation-finite. Since $\{\Ok{i}(C_E)\ |\
E\in\mc{E}, i\ge n\}\subseteq \bigcup_{i\ge n}\Ok{i}(\mc{C})$ and
the latter is representation-finite whenever $n\ge m$, we obtain
that $\bigcup_{i\ge n+1}\Ok{i}(\mc{E})$ is representation-finite
for $n=max\{m,k\}$. It follows that $\mc{E}$ is syzygy-finite.
%
%(2) If $\mc{C}$ is syzygy-finite from some $c$, then $\{\Ok{i}(C)\
%|\ i\le c, C\in\mc{C}\}$ is representation-finite for any
%$C\in\mc{C}$. By the proof of (1), we see that $\{\Ok{i}(E)\ |\
%i\le e, E\in\mc{E}\}$, where $e=\mr{min}\{c,a\}-1$, is
%representation-finite, for any $E\in\mc{E}$. Hence, $\mc{E}$ is
%syzygy-finite from $e$.
%

\ \hfill$\Box$

%
%%%%%%%%%%%%%%%%%%%%%%%%%%%%%%%%%%%%%%%%%%%%%%%%%%%%%%%%%%%%%%%%%%%%%
%%
%\vskip 30pt
%%
%\section{Syzygies and derived equivalences}
%

%%%%%%%%%%%%%%%%%%%%%%%%%%%%%%%%%%%%%%%%%%%%%%%%%%%%%
\zskip\

We need the following result on basic properties of tilting
complexes.

\bg{Lem}\label{t}%
Assume that there is an equivalence $\mc{F}:\ \dbr
\rightleftarrows\ \dbs\ :\mc{G}$. Let $T:=\mc{G}(S)$. Assume that
$T\in\kr{[a,k]}$ for some integers $a\le k$. Then

$(1)$ $\mc{F}(\dr{[c,d]})\subseteq\ds{[a-d,k-c]}$, for any $c\le
d$.
%If $M\in\kfr$ has the homological interval $[c,d]$. Then
%$\mc{F}(M)$ has the homological interval $[a-(d-c+1),b+(d-c+1)]$.

$(2)$  $\mc{G}(\ks{\mr{r}[c,d]})\subseteq\kr{\mr{r}[a+c,k+d]}$,
for any $c\le d$.
%If $B\in\kbr$ has the representation interval $[c,d]$. Then
%$\mc{F}(B)$ has the representation interval $[a+c,b+d]$.
%
\ed{Lem}

\Pf. (1) For any $M\in\dbr$ and any $i$, we have isomorphisms:

\sskip\

\hskip 109pt $\mr{H}^i(\mc{F}(M))\simeq \Hs(S,\mc{F}(M)[i])$

\hskip 166pt $\simeq\Hr(\mc{G}(S),\mc{G}(\mc{F}(M))[i])$

\hskip 166pt $\simeq \Hr(T,M[i]).$

\sskip\

Since $T\in\kr{[a,k]}$ and $M\in\dr{[c,d]}$, we obtain that
$\Hr(T,M[i])=0$ for $i\not\in [a-d,k-c]$. %[a-(d-c+1), b+(d-c+1)]$.
The conclusion then follows.

(2) Take any $B\in\ks{\mr{r}[c,d]}$. Note that there are triangles
$\Ok{i+1}(B)\to P_{i}\to \Ok{i}(B)\to$ for all integers $i$, where
each $P_{i}$ is a projective $S$-module, by Lemma \ref{Stri}.
Since $B\in\ks{\mr{r}[c,d]}$, we see that $\Ok{d}(B)$ is
(isomorphic to) a projective $S$-module and that
$B\simeq\Ok{c}(B)[c]$ by Lemma \ref{So} (2).
%$\sigma_{(\infty,c]}(B)=B$.
Applying the functor $\mc{G}$ to these triangles, we obtain
triangles $\mc{G}(\Ok{i+1}(B))\to \mc{G}(P_{i})\to
\mc{G}(\Ok{i}(B))\to$. Note that $\mc{G}(\Ok{d}(B)),
\mc{G}(P_{i})\in\mr{add}(\mc{G}(S))=\mr{add}T\subseteq
\kr{\mr{r}[a,k]}$, so we obtain that $\mc{G}(\Ok{c}(B))\in
\kr{\mr{r}[a,k+(d-c)]}$ from the above triangles, by using the
construction of cones.
%Since $\Ok{c}(B)=(\sigma_{(c,\infty]}(B))[-c]=B[-c]$,
Consequently, we see that
$\mc{G}(B)\simeq\mc{G}((\Ok{c}(B))[c])\simeq(\mc{G}(\Ok{c}(B)))[c]\in
\kr{\mr{r}[a+c,k+d]}$.
\ \hfill$\Box$

%%%%%%%%%%%%%%%%%%%%%%%%%%%%%%%%%%%%%%%%%%%%%%%%%%%%
\zskip\

The following result shows that derived equivalences preserve
syzygy-finite classes.

\bg{Pro}\label{syzc}%
Assume that there is an equivalence $\mc{F}: \dbr
\rightleftarrows\ \dbs\ :\mc{G}$. Let
$T:=\mc{G}(S)\in\kr{\mr{r}[a,k]}$ and  $\mc{C}\subseteq \dbs$.
 If $\mc{C}$  is
syzygy-finite, then $\mc{G(C)}$ is also syzygy-finite.
%
%$(2)$ If $\mc{C}$ is syzygy-finite from $c$, then $\mc{F(C)}$ is
%syzygy-finite from $e_{a,c}$, where $e_{a,c}=a$.
%
\ed{Pro}

\Pf.  Since $\bigcup_{i\ge n}\Ok{i}(\mc{C})$ is representation
finite for some $n$, we have some $M\in\dbs$ such that
$\bigcup_{i\ge n}\Ok{i}(\mc{C})\subseteq\mr{add}M$.
% Then
%$\Ok{i}(M)\in\mr{add}(M)$ for all $i$ and $M=\Ok{0}(M)$.

\bg{verse} $Claim$: $\mc{G}(\bigcup_{i\ge n}\Ok{i}(\mc{C}))$ is
$k$-syzygy-finite.

{\it Proof}.\ \ Take any $C\in\mc{C}$. Note that, for any $i$,
there is a triangle $\Ok{i+1}(C)\to C_{i}\to\Ok{i}(C)\to$, where
$C_{i}$ is projective, by Lemma \ref{Stri}. So, by applying the
functor $\mc{G}$, we obtain a triangle $\mc{G}(\Ok{i+1}(C))\to
\mc{G}(C_{i})\to\mc{G}(\Ok{i}(C))\to$. Note that
$\mc{G}(C_{i})\in\mr{add}T\subseteq\kr{\mr{r}[a,k]}$, so we have
that $\Ok{j}(\mc{G}(\Ok{i+1}(C))\simeq \Ok{j+1}(\mc{G}(\Ok{i}(C))$
for all $j\ge k$, by Lemma \ref{syzl}. It follows that

\hskip 50pt $\Ok{j}(\mc{G}(\Ok{i}(C)))\simeq
\Ok{j-1}(\mc{G}(\Ok{i+1}(C)))\simeq\cdots \simeq
\Ok{k}(\mc{G}(\Ok{j+i-k}(C)))$

\noindent for all $j\ge k$. Hence, we obtain that the class
$\bigcup_{j\ge k}\Ok{j}(\mc{G}(\bigcup_{i\ge
n}\Ok{i}(\mc{C})))=\bigcup_{j\ge k}\bigcup_{i\ge
n}\Ok{j}(\mc{G}(\Ok{i}(\mc{C})))$ is representation finite if and
only if the class $\bigcup_{j\ge k}\bigcup_{i\ge
n}\Ok{k}(\mc{G}(\Ok{j+i-k}(\mc{C})))=\bigcup_{j\ge
b}\Ok{k}(\mc{G}(\bigcup_{i\ge n}\Ok{j+i-k}(\mc{C})))$ is
representation finite. Since $\bigcup_{i\ge
n}\Ok{i}(\mc{C})\subseteq\mr{add}M$, we see that the last class is
contained in the class $\bigcup_{j\ge
b}\Ok{k}(\mc{G}(\mr{add}M))=\Ok{k}(\mc{G}(\mr{add}M))$ and hence
is representation finite. The claim then follows. \ed{verse}

Now take any $C\in\mc{C}$. Note that there is some $m_C$ such that
$C[-m_C]\simeq\Ok{m_C}(C)$ by Lemma \ref{So} (2),
%
% and consider it in $\kfs$, then there
%is some $m_C$ such that $C=\sigma_{[m_C,\infty)}(C)$.
%
so we have a triangle
 $(\Ok{n}(C))[n-1]\to B\to C\to $ by Lemma \ref{Stri}, where
$B\in\ks{[m_C,n-1]}$ and $n\ge m_C$. Then we obtain a triangle
$\mc{G}((\Ok{n}(C))[n-1])\to \mc{G}(B)\to\mc{G}(C)\to$, by
applying the functor $\mc{G}$. By Lemma \ref{t} (2),
$\mc{G}(B)\in\kr{[a+m_C,k+n-1]}\subseteq\kr{(-\infty,k+n-1]}$.
Note that $\mc{G}(\Ok{n}(\mc{C}))$ is $k$-syzygy-finite followed
from the claim, so
$\mc{G}((\Ok{n}(C))[n-1])=\mc{G}(\Ok{n}(C))[n-1]$ is
$(k+n-1)$-syzygy-finite, by Lemma \ref{So} (3). Hence, by
Proposition \ref{syzbj}, we have that $\mc{G}(C)$ is
$(k+n)$-syzygy-finite.
\ \hfill$\Box$

%%%%%%%%%%%%%%%%%%%%%%%%%%%%%%%%%%%%%%%%%%%%%%%%%%%%
\zskip\

Now we prove that derived equivalences preserve syzygy-finite
algebras.

\bg{Th}\label{Sfde}%Syz-finite-deriv-equi
Assume that $R,S$ are derived equivalent algebras. If $S$ is
syzygy-finite, then $R$ is also syzygy-finite.
\ed{Th}

\Pf.  By assumption, there is an equivalence $\mc{F}: \dbr
\rightleftarrows\ \dbs\ :\mc{G}$. Assume that
$T:=\mc{G}(S)\in\kr{\mr{r}[a,k]}$. Then
$\mc{F}(\mr{mod}R)=\mc{F}(\dr{[0,0]})\subseteq \ds{[a,k]}$, by
Lemma \ref{t} (1). If $S$ is syzygy-finite, then $\ds{[a,k]}$ is
syzygy-finite by Theorem \ref{Sfd}. It follows that
$\mc{F}(\mr{mod}R)$ is syzygy-finite. Hence, we obtain that
$\mr{mod}R=\mc{G(F}(\mr{mod}R))$ is syzygy-finite, by Proposition
\ref{syzc}. Thus, $R$ is syzygy-finite.
\ \hfill$\Box$

%%%%%%%%%%%%%%%%%%%%%%%%%%%%%%%%%%%%%%%%%%%%%%%%%%%%
\zskip\

For example, if $R$ is derived equivalent to a minimal
representation-infinite algebra or a monomial algebra, then $R$ is
syzygy-finite by the above theorem. In particular, $R$ is
CM-finite in the case.

Note that derived equivalences preserve Gorenstein algebras and
that syzygy-finiteness coincides with CM-finiteness for Gorenstein
algebras, we obtain the following corollary.

\bg{Cor}\label{CM}%Syz-finite-deriv-equi
Assume that $R,S$ are derived equivalent algebras. If $S$ is
Gorenstein CM-finite, then $R$ is also Gorenstein CM-finite.
\ed{Cor}

%%%%%%%%%%%%%%%%%%%%%%%%%%%%%%%%%%%%%%%%

\zskip\

\section{Igusa-Todorov algebras}

\hskip 20pt

%Now we turn to Igusa-Todorov algebras. Let $\mc{C}\subseteq \dbr$.
%We say that $\mc{C}$ is $n$-Igusa-Todorov provided that there is
%an $A$-module $_AV$ and some nonnegative integer $n$ such that,
%for any $M\in \Omega^n\mc{C}$, there is an exact sequence $0\to
%V_1\to V_0\to M\to 0$ with $V_0, V_1\in \mr{add}_AV$. The module
%$_AV$ is then called an $(n$-$)\mc{C}$-Igusa-Todorov module.Recall
%that

%Now let us recall the definition of Igusa-Todorov algebras
%introduced in [\ref{Wit}]. By definition

Recall that  an artin algebra $R$ is called $n$-{\it
Igusa-Todorov} provided that there exists a fixed $R$-module $V$
and a nonnegative integer $n$ such that, for any $M\in \mr{mod}R$,
there is an exact sequence $0\to V_1\to V_0\to \Ok{n}M\to 0$ with
$V_0, V_1\in \mr{add}V$ [\ref{Wit}].
% The class of Igusa-Todorov algebras is
%large and includes all algebras of representation dimension at
%most 3, algebras with radical cube zero and syzygy-finite algebras
%etc. [\ref{Wit}].
 A remarkable property of Igusa-Todorov algebras
that they satisfy the finitistic dimension conjecture. %In fact,
%most algebras which were recently proved to satisfy the finitistic
%dimension conjecture are Igusa-Todorov algebras.

 The following algebras are Igusa-Todorov.

\zskip\

$\bullet$ Algebras with radical cube zero.

$\bullet$ Algebras with representation dimension at most three.

$\bullet$ Syzygy-finite algebras.

$\bullet$ Algebras which are endomorphism algebras of modules over
representation-finite algebras [\ref{Wit}].

$\bullet$ Algebras with an ideal $I$ of finite projective dimensin
such that $I\mr{rad}^2R=0$ (or $I^2\mr{rad}R=0$) and $R/I$ is
syzygy-finite [\ref{Ws}].

%$\bullet$ Algebras with an ideal $I$ of finite projective dimensin
%such that $I\mr{rad}^2R=0$ (or $I^2\mr{rad}R=0$) and $R/I$ is
%syzygy-finite [\ref{Ws}].

$\bullet$ Algebras possessing a left idealized extension which is
$2$-syzygy-finite [\ref{Wit}].

$\bullet$ Algebras which are endomorphism algebras of projective
modules over $2$-Igusa-Todorov algebras  [\ref{Wit}].

$\bullet$ Algebras with an ideal $I$ of finite projective dimensin
such that $I\mr{rad}R=0$ and $R/I$ is Igusa-Todorov [\ref{Ws}].

\zskip\

Let $\mc{C}$ be a subclass of $\dbr$. We say that $\mc{C}$ is
relative hereditary provided that there is a complex $V\in\dbr$
such that, for any $M\in \mc{C}$, there is a triangle $V_1\to
V_0\to M\to$ with $V_0, V_1\in \mr{add}V$. We say that $\mc{C}$ is
an $n$-Igusa-Todorov class, for some integer $n$, provided that
$\Ok{n}(\mc{C})$  is relative hereditary. It is easy to see that
$\mc{C}$ is an Igusa-Todorov (resp., relative hereditary) class if
and only if $\mc{C}[n]$ is an Igusa-Todorov (resp., relative
hereditary)  class for any/some $n$.  It is also obvious that if
$\mc{E}\subseteq\mc{C}$ and $\mc{C}$ is Igusa-Todorov (resp.,
relative hereditary) then $\mc{E}$ is also Igusa-Todorov (resp.,
relative hereditary).

%%%%%%%%%%%%%%%%%%%%%%%%%%%%%%%%%%%%%%%%%%%%

\bg{Lem}\label{ITc}%
Let $\mc{C}\subseteq \dbr$. Then $\mc{C}$ is an Igusa-Todorov
class if and only if $\Ok{n}(\mc{C})$ is an Igusa-Todorov class
for any/some $n$.
\ed{Lem}

\Pf. We first prove that if $\Ok{m}(\mc{C})$ is relative
hereditary, then $\Ok{n}(\mc{C})$ is also relative hereditary, for
any $n\ge m$. In fact, by definition there is a complex $V\in\dbr$
such that, for any $M\in \Ok{m}(\mc{C})$, there is a triangle
$V_1\to V_0\to M\to$ with $V_0, V_1\in \mr{add}V$. By Proposition
\ref{Stn}, we have a triangle $\Ok{n-m}(V_1)\to \Ok{n-m}(V_0)\to
\Ok{n-m}(M)\to$. Note that
$\Ok{n-m}(V_1),\Ok{n-m}(V_0)\in\Ok{n-m}(\mr{add}V)$ is independent
of $M$, and $\Ok{n}(\mc{C})=\Ok{n-m}(\Ok{m}(\mc{C}))$ for $n\ge m$
by Lemma \ref{So} (4), so we obtain that $\Ok{n}(\mc{C})$ is
relative hereditary, for any $n\ge m$, by definition.

Now assume that $\mc{C}$ is an $m$-Igusa-Todorov algebra. Then
$\Ok{m}(\mc{C})$ is relative hereditary. Let $n$ be an integer,
and take an integer $t>\mr{max}\{0,m-n\}$, then we have that
$\Ok{t}(\Ok{n}(\mc{C}))=\Ok{n+t}(\mc{C})$ is relative hereditary,
by the above argument. It follows that $\Ok{n}(\mc{C})$ is an
Igusa-Todorov algebra by definition. Conversely, assume that
$\Ok{n}(\mc{C})$ is a $t$-Igusa-Todorov algebra, for some integers
$n,t$, then $\Ok{n+t}(\mc{C})$ is relative hereditary. Hence
$\mc{C}$ is an Igusa-Todorov algebra by definition.
\ \hfill$\Box$

%%%%%%%%%%%%%%%%%%%%%%%%%%%%%%%%%%%%%%%%%%%%
\zskip\

The following result gives a characterization of Igusa-Todorov
algebras in term of derived categories.

\bg{Th}\label{ITd}%
An algebra $R$ is Igusa-Todorov if and only if $\dr{[a,k]}$ is an
Igusa-Todorov class for any/some $a\le k$.
\ed{Th}

\Pf. If $R$ is an Igusa-Todorov algebra, then $\mr{mod}R$ is an
Igusa-Todorov class by definition. Note that
$\mr{mod}R=\Ok{k}(\dr{[a,k]})$, so we further have that
$\dr{[a,k]}$ is an Igusa-Todorov class, by Lemma \ref{ITc}.

Conversely, assume that $\dr{[a,k]}$ is an Igusa-Todorov class,
then $\mr{mod}R$ is also an Igusa-Todorov class by Lemma \ref{ITc}
again. Thus there is a fixed complex $V\in\dbr$ and an integer $v$
such that, for any $M\in\mr{mod}R$, there is a triangle $V_1\to
V_0\to \Ok{v}(M)\to$ with $V_1,V_0\in\mr{add}V$. Assume that
$V\in\dr{[c,d]}$ for some $c\le d$, then $\Ok{n}(V)\in\mr{mod}R$
for all $n\ge \mr{max}\{d,0\}$, by Lemma \ref{So} (2). By
Proposition \ref{Stn} and Lemma \ref{So} (4), there is a triangle
$\Ok{n}(V_1)\to \Ok{n}(V_0)\to \Ok{n}(\Ok{v}(M))
(\simeq\Ok{n+v}(M))\to$. Take the integer $n$ such that $n\ge
\mr{max}\{0,d,-v\}$, then all terms in the last triangle are
$R$-modules. Hence, we have an exact sequence $0\to \Ok{n}(V_1)\to
\Ok{n}(V_0)\to \Ok{n+v}(M)\to 0$. Note that
$\Ok{n}(V_1),\Ok{n}(V_0)\in\Ok{n}(\mr{add}V)\subseteq \mr{mod}R$
and $n+v$ are independent of $M$, so we obtain that  $R$ is an
Igusa-Todorov algebra by definition.
\ \hfill$\Box$

%%%%%%%%%%%%%%%%%%%%%%%%%%%%%%%%%%%%%%%%%%%%%%%%%%%%
\zskip\

The following result shows that derived equivalences preserve
Igusa-Todorov classes.

\bg{Pro}\label{ITcd}%
Assume that there is an equivalence $\mc{F}: \dbr
\rightleftarrows\ \dbs\ :\mc{G}$. Let
$T:=\mc{G}(S)\in\kr{\mr{r}[a,k]}$ and  $\mc{C}\subseteq \dbs$.
 If $\mc{C}$ is an Igusa-Todorov
class, then $\mc{G(C)}$ is also an Igusa-Todorov class.
\ed{Pro}

\Pf. Assume that $\mc{C}$ is an $n$-Igusa-Todorov class for some
$n$, then there is a fixed complex $V\in\dbr$ such that, for any
$C\in\mc{C}$, there is a triangle $V_1\to V_0\to \Ok{n}(C)\to$
with $V_1,V_0\in\mr{add}V$. Hence, we also have a triangle
$\mc{G}(V_1)\to \mc{G}(V_0)\to \mc{G}(\Ok{n}(C))\to$, by applying
the functor $\mc{G}$.
 Similarly as in the proof of Proposition
\ref{syzc}, there is a triangle $\mc{G}((\Ok{n}(C))[n-1])\to
\mc{G}(B)\to \mc{G}(C)\to$ with
$\mc{G}(B)\in\kr{\mr{r}(-\infty,k+n-1]}$. So, we have that

\sskip\ \noindent \hskip 80pt
$\Ok{k+n}(\mc{G}(C))\simeq\Ok{k+n-1}(\mc{G}((\Ok{n}(C))[n-1]))$

\sskip\
 \noindent \hskip 144pt$\simeq \Ok{k+n-1}(\mc{G}((\Ok{n}(C)))[n-1])\simeq
\Ok{k}(\mc{G}((\Ok{n}(C)))$,

\sskip\

\noindent by Proposition \ref{syzl} and Lemma \ref{So} (3). Note
that we have a triangle $\Ok{k}(\mc{G}(V_1))\to
\Ok{k}(\mc{G}(V_0))\to \Ok{k+n}(\mc{G}(C))\to$, by Proposition
\ref{Stn} and the above argument. Since $b+n$ and
$\Ok{k}(\mc{G}(V_1)),
\Ok{k}(\mc{G}(V_0))\in\Ok{k}(\mc{G}(\mr{add}V))$ are independent
of $C$, so we obtain that $\Ok{k+n}(\mc{G}(C))$ is relative
hereditary, i.e., $\mc{G(C)}$ is a $(k+n)$-Igusa-Todorov class.
\ \hfill$\Box$

%%%%%%%%%%%%%%%%%%%%%%%%%%%%%%%%%%%%%%%%%%%%
\zskip\

Now we obtain the following important result on Igusa-Todorov
algebras.

\bg{Th}\label{ITde}%
Assume that $R,S$ are derived equivalent algebras. If $S$ is an
Igusa-Todorov algebra, then $R$ is also an Igusa-Todorov algebra.
\ed{Th}

\Pf. The proof is similar to that of Theorem \ref{Sfde}. Namely,
if the equivalence is given by $\mc{F}: \dbr \rightleftarrows\
\dbs\ :\mc{G}$, then we can obtain that $\mc{F}(\mr{mod}R)$ is an
Igusa-Todorov class provided that $S$ is an Igusa-Todorov algebra.
Hence, we obtain that $\mr{mod}R=\mc{G(F}(\mr{mod}R))$ is  an
Igusa-Todorov class, by Proposition \ref{ITcd}. Thus, $R$ is  an
Igusa-Todorov algebra by Theorem \ref{ITd}.
\ \hfill$\Box$

%%%%%%%%%%%%%%%%%%%%%%%%%%%%%%%%%%%%%%%%%%%%%%%%%%%%%%%%%%%%%%%%%%%%%
\zskip\

For instance, if $R$ is derived equivalent to an algebra with
radical cube zero, or an algebra with representation dimension at
most three, then $R$ must be an Igusa-Todorov algebra.

\bskip\

\bskip\

\section{Auslander's condition}

\hskip 15pt

Let $M\in\dbr$ and $\mc{C}\subseteq \dbr$. We define the
$\mc{C}$-Auslander bound of $M$ to be the minimal integer $m$, or
$\infty$ if such minimal integer doesn't exist, such that
$\Hr(M,N[i])=0$ for all $i>m$, whenever $N\in\mc{C}$ satisfies
that $\Hr(M,N[i])=0$ for all but finitely many $i$. Let
$\mc{E}\subseteq \dbr$. The global $\mc{C}$-Auslander bound of
$\mc{E}$ is the supremum of all $\mc{C}$-Auslander bounds of
objects in $\mc{E}$. The finitistic $\mc{C}$-Auslander bound of
$\mc{E}$ is the supremum of  all $\mc{C}$-Auslander bounds of
objects in $\mc{E}$ whose  $\mc{C}$-Auslander bound is finite.

In case that $M$ is an $R$-module and $\mc{C}=\mc{E}=\mr{mod}R$,
the notions given here coincide with the usual ones in module
categories [\ref{CH}, \ref{Wab}], i.e., Auslander bound of $M$,
global Auslander bound of the algebra $R$ and the finitistic
Auslander bound of the algebra $R$, respectively.

We have the following easy observation.

\bg{Lem}\label{ACl}%
Let $M\in\dbr$ and $\mc{C}\subseteq \dbr$. Then $M$ has finite
$\mc{C}$-Auslander bound if and only if $M[m]$ has finite
$(\mc{C}[n])$-Auslander bound for any/some integers $m,n$.
\ed{Lem}

It is also easy to see that if a complex $M$ has finite
$\mc{C}$-Auslander bound, then $M$ also has finite
$\mc{E}$-Auslander bound for any $\mc{E}\subseteq \mc{C}$.

An algebra $R$ is called an AC-algebra provided that every
$R$-module has finite $(\mr{mod}R)$-Auslander bound [\ref{CH}].
Auslander conjecture asserts that all algebras are AC-algebras.
However, the conjecture fails in general [\ref{JS},\ref{Sm}]. We
refer to [\ref{CH}] for the list of AC-algebras. In [\ref{Wab}],
the author suggests a revisited version of Auslander conjecture,
named the finitistic Auslander conjecture, which asserts that the
finitistic Auslander bound of every algebra is finite. Note that
the finitistic Auslander conjecture implies the finitistic
dimension conjecture.

We have the following characterization of AC-algebras in term of
derived categories.

\bg{Th}\label{ACd}%AC-deriv
$(1)$ $R$ is an AC-algebra if and only if every complex in $\dbr$
has finite $(\dr{[c,d]})$-Auslander bound for any/some integers
$c\le d$.

$(2)$ The global Auslander bound of $R$ is finite if and only if,
for any/some integers $a\le k$ and $c\le d$, the global
$(\dr{[c,d]})$-Auslander bound of  $\dr{[a,k]}$ is finite.

$(3)$ The finitistic Auslander bound of $R$ is finite if and only
if, for any/some integers $a\le k$ and $c\le d$, the finitistic
$(\dr{[c,d]})$-Auslander bound of  $\dr{[a,k]}$ is finite.
\ed{Th}

\Pf. (1) The if-part. Assume that every complex in $\dbr$ has
finite $(\dr{[c,d]})$-Auslander bound for some integers $c\le d$.
Then every complex has finite $((\mr{mod}R)[c])$-Auslander bound
since $(\mr{mod}R)[c]\subseteq\dr{[c,d]}$. It follows that every
complex, in particular every $R$-module, has finite
$(\mr{mod}R)$-Auslander bound by Lemma \ref{ACl}. Hence, $R$ is an
AC-algebra.

The only-if part. Take any $M\in\dbr$ and any integers $c\le d$..
Assume that $M\in\dr{[a,k]}$ for some integers $a\le k$, then we
have that $\Ok{k}(M)$ is an $R$-module and $M\simeq\Ok{a}(M)[a]$,
by Lemma \ref{So} (2). Now take any $N\in\dr{[c,d]}$, then we
obtain that

\sskip\

\hskip 100pt $\Hr(M,N[j])\simeq\Hr(\Ok{a}(M)[a],N[j])$

\sskip\

\hskip 50pt $\simeq\Hr(\Ok{a}(M),N[j-a])\stackrel{j-k>-c}{\simeq}
\Hr(\Ok{k}(M),N[j-k])$

\sskip\

\noindent  for all $j>k-c$, by dimension shifting in Proposition
\ref{St}. Now note that $\Ok{k}(M)$ is an $R$-module
%the cosyzygy $\Oki{c}(N)$ is an $R$-module
and $N=\Oki{d}(N)[d]$ by Lemma \ref{So}' (2), so we further obtain
that

\sskip\

\hskip 50pt
$\Hr(\Ok{k}(M),N[j-k])\simeq\Hr(\Ok{k}(M),(\Oki{d}(N))[j-k+d]) $

\sskip\

\hskip 70pt $\stackrel{j-k+c>0}{\simeq}
\Hr(\Ok{k}(M),(\Oki{c}(N))[j-k+c])$,

\sskip\

\noindent for all $j>k-c$,  by dimension shifting in Proposition
\ref{St}'. It follows that

\sskip\

\hskip 80pt $\Hr(M,N[j])\simeq
\Hr(\Ok{k}(M),(\Oki{c}(N))[j-k+c])$,

\sskip\

\noindent for all $j>k-c$.%t:=\mr{max}\{b-c,b+2d-2a-c\}$.

Hence, if $\Hr(M,N[j])=0$ for all but finitely many $j$, then we
have that  $\Hr(\Ok{k}(M),  (\Oki{c}(N))[i])=0$ for all but
finitely many $i$. Since $R$ is an AC-algebra and
$\Ok{k}(M),\Oki{c}(N)\subseteq \mr{mod}R$, there is an integer
$l\ge 0$, independent of $\Ok{k}(N)$, such that
$\Hr(\Ok{k}(M),(\Oki{c}(N))[i])=0$ for all $i>l$. It follows from
the above isomorphism that $\Hr(M,N[j])=0$ for all $j$ such that
$j-k+c>l$ and $j>k-c$, i.e., for all $j>k-c+l$. Note that $k-c+l$
is independent of $N$, so we obtain that $M$ has finite
$(\dr{[c,d]})$-Auslander bound.

(2) and (3). The proofs are similar to (1), just note that
Auslander bounds are unique in both cases.
\ \hfill$\Box$

%%%%%%%%%%%%%%%%%%%%%%%%%%%%%%%%%%%%%%%%%%%%%%%%%%%%%
\zskip\

Then, we obtain an important property of AC-algebras and algebras
satisfying the finitistic Auslander conjecture.

\bg{Th}\label{ACde}%
Assume that $R,S$ are derived equivalent algebras.

$(1)$ If $S$ is an AC algebra, then $R$ is also an AC algebra.

$(2)$ If the global Auslander bound of $S$ is finite, then  the
global Auslander bound of $R$ is also finite.

$(3)$ If the finitistic Auslander bound of $S$ is finite, then the
finitistic Auslander bound of $R$ is also finite.
\ed{Th}

\Pf. (1) Take any $M,N\in\mr{mod}R$. If $\Hr(M,N[i])=0$ for all
but finitely many $i$, then $\Hs(\mc{F}(M),\mc{F}(N)[i])=0$ for
all but finitely many $i$. Note that
$\mc{F}(N)\in\mc{F}(\mr{mod}S)\subseteq \ds{[a,k]}$ for some fixed
integer $a\le k$, by Lemma \ref{t}. Since  $S$ is an AC algebra,
we have that $\mc{F}(M)\in\dbs$ has finite
$(\ds{[a,k]})$-Auslander bound by Theorem \ref{ACde}. Then there
is an integer $m$, independent of $N$, such that
$\Hs(\mc{F}(M),\mc{F}(N)[i])=0$ for all $i>m$. It follows that
$\Hr(M,N[i])=\Hr(\mc{GF}(M),\mc{GF}(N)[i])=0$ for all $i>m$, i.e.,
$M$ has finite $(\mr{mod}R)$-Auslander bound. Hence, $R$ is an
AC-algebra.

(2) and (3). The proofs are similar.
\ \hfill$\Box$

%%%%%%%%%%%%%%%%%%%%%%%%%%%%%%%%%%%%%%%%%%%%%%%%%%%%%
\zskip\

Now we turn to algebras satisfying the generalized
Auslander-Reiten conjecture. We note that an equivalent statement
of the generalized Auslander-Reiten conjecture is  that if $M$ is
an $R$-module such that $\Hr(M,(M\oplus R)[i])=0$ for all but
finitely many $i$, then $M$ is of finite projective dimension,
i.e., $M$ is (isomorphic to) a complex in $\kbr$.

%, see [\ref{Wab}]. This
%statement is equivalent to that in the introduction.

\bg{Lem}\label{Gcl}%
The following conditions are equivalent for a complex $M\in\dbr$.

$(1)$ $\Hr(M,(M\oplus R)[i])=0$ for all but finitely many $i$.

$(2)$ $\Hr(\Ok{n}(M),(\Ok{n}(M)\oplus R)[i])=0$ for any/some
integer $n$ and all but finitely many $i$.

$(3)$ $\Hr(M,(M\oplus T)[i])=0$ for any/some tilting complex $T$
and all but finitely many $i$.
\ed{Lem}

\Pf. (1) $\Leftrightarrow$ (2) Note that $M=\Ok{m}(M)[m]$ for some
integer $m$, so we obtain that $\Hr(\Ok{m}(M),(\Ok{m}(M)\oplus
R)[i])=0$ for all but finitely many $i$ if and only if $\Hr(M$,
$(M\oplus R)[i])=0$ for all but finitely many $i$. By Proposition
\ref{St}, there is a triangle $\Ok{j+1}(M)\to P_j\to \Ok{j}(M)\to$
with $P_j$ projective, for
any $j$. So, one can easily check that % we have  that
$\Hr(\Ok{j}(M), (\Ok{j}(M)\oplus R)[i])=0$ for all but finitely
many $i$ if and only if $\Hr(\Ok{j+1}(M),(\Ok{j+1}(M)\oplus
R)[i])=0$ for all but finitely many $i$. The conclusion then
follows.

(1) $\Leftrightarrow$ (3) If $T$ is a tilting complex, then $R$
generates $T$ and $T$ generates $R$ in $\kbr$, both via finitely
many steps. Also note that $\Hr(B,M[i])=0$ for all but finitely
many $i$, whenever $B\in\kbr$. It follows that $\Hr(M,(M\oplus
T)[i])=0$ for all but finitely many $i$ if and only if
$\Hr(M,(M\oplus R)[i])=0$ for  all but finitely many $i$.
\ \hfill$\Box$

%%%%%%%%%%%%%%%%%%%%%%%%%%%%%%%%%%%%%%%%%%%%%%%%%%%%%%%%%%
\zskip\

Now we can provide a derived version of the generalized
Auslander-Reiten conjecture.

\bg{Th}\label{Gcd}%
An algebra $R$ satisfies the generalized Auslander-Reiten
conjecture if and only if, for any $M\in\dbr$ such that
$\Hr(M,(M\oplus T)[i])=0$ for any/some tilting complex $T$ and all
but finitely many $i$, it holds that $M\in\kbr$.
\ed{Th}

\Pf. The sufficient part follows from Lemma \ref{Gcl}.

The necessary part. Note that there is some $n$ such that
$\Ok{n}(M)$ is an $R$-module, by Lemma \ref{So} (2). If $M$
satisfies that $\Hr(M,(M\oplus R)[i])=0$ for some tilting complex
and all but finitely many $i$, then we have that
$\Hr(\Ok{n}(M),(\Ok{n}(M)\oplus R)[i])=0$ for all but finitely
many $i$, by Lemma \ref{Gcl}. Since $R$  satisfies the generalized
Auslander-Reiten conjecture, we obtain the $\Ok{n}(M)\in\kbr$. It
follows that $M\in\kbr$, by Lemma \ref{So} (6).
 \ \hfill$\Box$

%%%%%%%%%%%%%%%%%%%%%%%%%%%%%%%%%%%%%%%%%%%%%%%%%%%%%%%%%%
\zskip\

Then, we can show that derived equivalences preserve generalized
Auslander-Reiten conjecture.

\bg{Th}\label{Gcde}%
Assume that $R,S$ are derived equivalent algebras. If $S$
satisfies the generalized Auslander-Reiten conjecture, then $R$
also satisfies the generalized Auslander-Reiten conjecture.
\ed{Th}

\Pf. Take any $M\in\dbr$ such that $\Hr(M,(M\oplus R)[i])=0$ for
all but finitely many $i$. Then $\Hs(\mc{F}(M),(\mc{F}(M)\oplus
\mc{F}(R))[i])=0$ for all but finitely many $i$. Note that
$\mc{F}(R)$ is a tilting complex in $\dbs$ and that $S$ satisfies
the generalized Auslander-Reiten conjecture, so we obtain that
$\mc{F}(M)\in\kbs$ by Theorem \ref{Gcd}. It follows that
$M\simeq\mc{GF}(M)\in\kbr$. Hence $R$ satisfies the generalized
Auslander-Reiten conjecture by Theorem \ref{Gcd} again.
\ \hfill$\Box$

%%%%%%%%%%%%%%%%%%%%%%%%%%%%%%%%%%%%%%%%%%%%
%\bskip\
%
%\section{Examples}
%

%1. See "properly stratified algebras and tilting"

%%-------------------------------------------------------------------
%%%%%%%%%%%%%%%%%%%%%%%%%%%%%%%%%%%%%%%%%%%%%%%%%%%%%%%%%%%%%%%%%%%%%
%%%%
%%%%%%%%%%%%%%%%%%%%%%%%%%%%%%%%%%%%%%%%%%%%%%%%%%%%%%%%%%%%%%%%%%%%%%%%%%%%%%%%%%%%
%%%**************************ÖÂ Ð»**********************************************
%\vskip 30pt
%\begin{center}
%{\bf ACKNOWLEDGEMENTS}
%\end{center}
%\hskip 15pt The author thanks for the referee's highly valuable
%comments and suggestions. In particular, the notions of
%``syzygy-bounded" and ``syzygy-equivalent'' are due to the
%referee, as well as Propositions \ref{syeqhom}.
%%, which makes the
%%paper more readable and .

%%It is a pleasure to thank the referee for helpful suggestions and
%%%%%%%%%%%%%%%%%%%%%%%%%%%%%%%%%%%%%%%%%%%%%%%%%%%%%%%%%%%%%%%%%%%%%%%%%%%%%%%%%%
\vskip 30pt

{\small

}

\vskip 50pt{\it

%\noindent Address:

\noindent Jiaqun Wei

\noindent  School of Mathematics Science, Nanjing Normal
University, Nanjing 210046, China

\noindent Email: weijiaqun@njnu.edu.cn }

%
%\vskip30pt
%
%\noindent{\it Current address:}
%
%\vskip10pt \noindent
%Universitaet zu Koeln\\
%\noindent Mathematisches Institut\\
%\noindent Weyertal 86-90\\
%\noindent D-50931 Koeln\\
%\noindent Germany

\end{document}